\documentclass{agtart_a}
\pdfoutput=1

\usepackage[all]{xy}

\title[On submanifolds in locally symmetric spaces of noncompact type]{On submanifolds in locally symmetric spaces\\ of noncompact type}

\author{Jean-Francois Lafont}
\givenname{Jean-Francois}
\surname{Lafont}
\address{Department of Mathematics\\
The Ohio State University\\\newline
231 West 18th Avenue\\
Columbus, OH 43210-1174}
\email{jlafont@math.ohio-state.edu}
\urladdr{}

\author{Benjamin Schmidt}
\givenname{Benjamin}
\surname{Schmidt}
\address{Department of Mathematics\\
University of Chicago\\\newline
5734 S\,University Avenue\\
Chicago, IL 60637}
\email{schmidt@math.uchicago.edu}
\urladdr{}

\volumenumber{6}
\issuenumber{}  
\publicationyear{2006}
\papernumber{85}
\startpage{2455}
\endpage{2472}

\doi{}
\MR{}
\Zbl{}

\keyword{locally symmetric space}
\keyword{duality}
\keyword{tangential map}
\keyword{Matsushima's map}
\subject{primary}{msc2000}{53C35}
\subject{secondary}{msc2000}{57T15}
\subject{secondary}{msc2000}{55R37}
\subject{secondary}{msc2000}{57R42}
\subject{secondary}{msc2000}{57R45}

\received{10 July 2006}
\revised{}
\accepted{9 November 2006}
\published{15 December 2006}
\publishedonline{15 December 2006}
\proposed{}
\seconded{}
\corresponding{}
\editor{}
\version{}

\arxivreference{math.GT/0607242}

\let\xysavmatrix\xymatrix
\def\xymatrix{\disablesubscriptcorrection\xysavmatrix}
\AtBeginDocument{\let\bar\wbar\let\tilde\wtilde\let\hat\what}
\makeop{Sing}
\makeop{rk}
\makeop{ker}

\newcommand{\p}{\prime}
\newcommand{\mC}{\mathbb C}
\newcommand{\mZ}{\mathbb Z}
\newcommand{\mR}{\mathbb R}

\newcommand{\mH}{\mathbb H}
\newcommand{\mO}{\mathbb O}
\newcommand{\bs}{\backslash}
\newcommand{\G}{\Gamma}
\newcommand{\gO}{\Omega}
\newcommand{\gL}{\Lambda}
\newcommand{\frakg}{\mathfrak g}
\newcommand{\ft}{\mathfrak t}
\newcommand{\fp}{\mathfrak p}
\newcommand{\SL}{\mathit{SL}}
\newcommand{\SU}{\mathit{SU}}

\makeatletter
\def\cnewtheorem#1[#2]#3{\newtheorem{#1}{#3}[section]
\expandafter\let\csname c@#1\endcsname\c@Lem}
\makeatother

\makeautorefname{Lem}{Lemma}
\cnewtheorem{Prop}[Lem]{Proposition}
\makeautorefname{Prop}{Proposition}
\cnewtheorem{Thm}[Lem]{Theorem}
\makeautorefname{Thm}{Theorem}
\cnewtheorem{Cor}[Lem]{Corollary}
\makeautorefname{Cor}{Corollary}

\theoremstyle{definition}
\newtheorem*{Def}{Definition}

\theoremstyle{remark}
\newtheorem*{Rmk}{Remark}

\begin{document}

\begin{asciiabstract}
Given a connected, compact, totally geodesic submanifold Y^m of
noncompact type inside a compact locally symmetric space of noncompact
type X^n, we provide a sufficient condition that ensures that [Y^m] is
nonzero in H_m(X^n; R); in low dimensions, our condition is also
necessary.  We provide conditions under which there exist a tangential
map of pairs from a finite cover (X-bar,Y-bar) to the nonnegatively
curved duals (X_u,Y_u).
\end{asciiabstract}

\begin{htmlabstract}
Given a connected, compact, totally geodesic submanifold Y<sup>m</sup>
of noncompact type inside a compact locally symmetric space of
noncompact type X<sup>n</sup>, we provide a sufficient condition that
ensures that [Y<sup>m</sup>] is nonzero in H<sub>m</sub>(X<sup>n</sup>;
<b>R</b>); in low dimensions, our condition is also necessary.  We
provide conditions under which there exist a tangential map of pairs
from a finite cover (X-bar,Y-bar) to the nonnegatively curved duals
(X<sub>u</sub>,Y<sub>u</sub>).
\end{htmlabstract}

\begin{abstract}
Given a connected, compact, totally geodesic submanifold $Y^m$ of
noncompact type inside a compact locally symmetric space of noncompact
type $X^n$, we provide a sufficient condition that ensures that
$[Y^m]\neq 0\in H_m(X^n; \mathbb{R})$; in low dimensions, our
condition is also necessary.  We provide conditions under which there
exist a tangential map of pairs from a finite cover $(\bar X,\bar Y)$
to the nonnegatively curved duals $(X_u,Y_u)$.
\end{abstract}

\maketitle

\section{Introduction}

In this paper, we propose to study totally geodesic submanifolds
inside locally symmetric spaces. Let us start by fixing some
notation: $(X^n,Y^m)$ will always refer to a pair of compact locally
symmetric spaces of noncompact type, with $Y^m\subset X^n$ a
totally geodesic submanifold.  The spaces $X^n,Y^m$ will be locally
modelled on $G/K$, $G^\p/K^\p$ respectively, where $G,G^\p$ are a
pair of semisimple Lie groups, and $K,K^\p$ are a pair of maximal
compact subgroups in the respective $G,G^\p$. Note that, since
$Y^m\subset X^n$ is totally geodesic, one can view $G^\p$ as a
subgroup of $G$, and hence one can take $K^\prime = K\cap G^\prime$.
We will denote by $X_u=G_u/K$, $Y_u=G_u^\p/K^\p$ the nonnegatively
curved dual symmetric spaces to the nonpositively curved spaces $G/K$,
$G^\p/K^\p$.

Note that for a pair $(X^n,Y^m)$, the submanifold $Y^m$ is always
{\it homotopically\/} nontrivial.  Indeed, the inclusion induces a
monomorphism on the level of fundamental groups.  A more subtle
question is whether the submanifold $Y^m$ is {\it homologically\/}
nontrivial, ie whether $[Y^m]\neq 0 \in H_m(X^n ;\mR)$ (or in
$H_m(X^n ;\mZ)$).  Our first result provides a criterion for
detecting when a totally geodesic submanifold $Y^m$ is homologically
nontrivial (over $\mathbb R$) in $X^n$.

\begin{Thm}\label{Theorem1.1}
Let $Y^m\hookrightarrow X^n$ be a compact totally geodesic
submanifold of noncompact type inside a compact locally symmetric
space of noncompact type, and denote by $\rho$ the map on cohomology
$H^m(X_u;\mR)\rightarrow H^m(Y_u;\mR)\simeq \mR$ induced by the
embedding $Y_u\hookrightarrow X_u$.  Then we have the following:
\begin{itemize}
\item If $[Y^m]=0 \in H_m(X^n; \mR)$ then the map $\rho$ is identically zero.
\item If $\rho$ is identically zero, and $m\leq m(\frakg)$, where $m(\frakg)$ is the Matsushima
constant corresponding to the Lie algebra $\frakg$ of the Lie group $G$, then
we have that $[Y^m]=0 \in H_m(X^n; \mR)$.
\end{itemize}
\end{Thm}

Our proof of this first result is an adaptation of an argument of
Matsushima \cite{Mat} and relies on the existence of certain
compatible maps (the Matsushima maps) from the real cohomology of
the pair of nonnegatively curved duals $(X_u,Y_u)$ to the real
cohomology of the nonpositively curved pair $(X^n,Y^m)$.  It is
reasonable to ask whether this map can be realized {\it
geometrically}.  Our second result, extending work of Okun
\cite{O1}, shows that this can sometimes be achieved {\it
rationally}:

\begin{Thm}\label{Theorem1.2}
Assume that $Y^m\hookrightarrow X^n$ is a totally geodesic embedding
of compact, locally symmetric spaces of noncompact type.
Furthermore, assume that the map $G_u^\p\hookrightarrow G_u$ induced
by the inclusion $Y\hookrightarrow X$ is a $\pi_i$--isomorphism,
for $i< m$, and a surjection on $\pi_m$. Then there exists a finite
cover $\bar X$ of $X^n$, and a connected lift $\bar Y \subset \bar
X$ of $Y^m$, with the property that there exists a {\it tangential
map} of pairs $(\bar X, \bar Y) \rightarrow (X_u, Y_u)$.  If in
addition we have $\rk(G_u)=\rk(K)$ and $\rk(G_u^\p)=\rk(K^\p)$, then the
respective tangential maps induce the Matsushima maps on cohomology.
\end{Thm}

Since the tangent bundle of the submanifold $Y^m$ Whitney sum with
the normal bundle of $Y^m$ in $X^n$ yields the restriction of the
tangent bundle of $X^n$ to the submanifold $Y^m$, this gives the
immediate:

\begin{Cor}\label{Corollary1.3}
Under the hypotheses of \fullref{Theorem1.2}, we have that the
pullback of the normal bundle of $Y_u$ in $X_u$ is stably equivalent
to the normal bundle of $\bar Y^m$ in $\bar X^n$.
\end{Cor}

In the previous corollary, we note that if $2m+1\leq n$, then these
two bundles are in fact isomorphic (see for instance Husemoller \cite[Chapter 8,
Theorem 1.5]{H}).

An example where the hypotheses of the \fullref{Theorem1.2} are satisfied
arises in the situation where $Y^m$, $X^n$ are real hyperbolic
manifolds. Specializing \fullref{Corollary1.3} to this situation, we
obtain:

\begin{Cor}\label{Corollary1.4}
Let $Y^m\hookrightarrow X^n$ be a totally geodesic embedding, where
$X^n$, $Y^m$ are compact hyperbolic manifolds, and assume that
$2m+1\leq n$. Then there exists a finite cover $\bar X$ of $X^n$,
and a connected lift $\bar Y$ of $Y^m$, with the property that the
normal bundle of $\bar Y$ in $\bar X$ is trivial.
\end{Cor}

While the hypotheses of \fullref{Theorem1.2} are fairly technical, we
point out that there exist several examples of inclusions
$Y^m\hookrightarrow X^n$ satisfying the hypotheses of the theorem.
The proof of \fullref{Corollary1.4}, as well as a discussion of some
further examples will be included at the end of Section 4.  Finally,
we will conclude the paper with various remarks and open questions
in Section 5.

\vspace{-4pt}
\paragraph{Acknowledgements}  
This research was partially conducted during the period when B\,Schmidt was employed by the Clay Mathematics Institute as a Liftoff Fellow. The research of J\,-F\,Lafont was partly supported by the National Science Foundation under grant DMS - 0606002.  The authors would like to thank the anonymous referee
for pointing out a simplification in our original proof of \fullref{Proposition5.2}.
\vspace{-4pt}

\section{Background}
\vspace{-4pt}

In this section, we provide some discussion of the statements of our theorems.  We also
introduce some of the ingredients that will be used in the proofs of our results.
\vspace{-4pt}

\subsection{Dual symmetric spaces}
\vspace{-4pt}

Let us start by recalling the definition of dual symmetric spaces:
\vspace{-4pt}

\begin{Def}
Given a symmetric space $G/K$ of noncompact type, we define the
{\it dual symmetric space\/} in the following manner. Let $G_\mC$
denote the complexification of the semisimple Lie group $G$, and
let $G_u$ denote the maximal compact subgroup in $G_\mC$.  Since $K$
is compact, under the natural inclusions $K\subset G\subset G_\mC$,
we can assume that $K\subset G_u$ (up to conjugation). The symmetric
space dual to $G/K$ is defined to be the symmetric space $G_u/K$.
By abuse of language, if $X=\Gamma \backslash G/K$ is a locally
symmetric space modelled on the symmetric space $G/K$, we will say
that $X$ and $G_u/K$ are dual spaces.
\end{Def}
\vspace{-4pt}

Now assume that $Y^m\hookrightarrow X^n$ is a totally geodesic
submanifold, where both $Y^m$, $X^n$ are locally symmetric spaces of
noncompact type.  Fixing a lift of $Y$, we have a totally geodesic
embedding of the universal covers:
$$G^\p/K^\p =\tilde Y\hookrightarrow \tilde X=G/K$$
Corresponding to this totally geodesic embedding, we get a natural commutative diagram:
$$\xymatrix{G^\p \ar[r] & G\\
K^\p \ar[r] \ar[u] & K \ar[u]}$$
which, after passing to the complexification, and descending to the maximal compacts,
yields a commutative diagram:
$$\xymatrix{G^\p_u \ar[r] & G_u\\
K^\p \ar[r] \ar[u] & K \ar[u]}$$
In particular, corresponding to the totally geodesic embedding $Y\hookrightarrow X$,
we see that there is a totally geodesic embedding of the dual symmetric spaces $G_u^\p/K^\p
\hookrightarrow G/K$.

\subsection{Classifying spaces}

For $G$ a continuous group let $EG$ denote a contractible space
which supports a free $G$--action. The quotient space, denoted $BG$,
is called a {\it classifying space\/} for principal $G$--bundles. This
terminology is justified by the fact that, for any topological space
$X$, there is a bijective correspondence between (1) isomorphism
classes of principal $G$--bundles over $X$, and (2) homotopy classes
of maps from $X$ to $BG$.  Note that the spaces $EG$ are only
defined up to $G$--equivariant homotopies, and likewise the spaces
$BG$ are only defined up to homotopy.  Milnor \cite{Mil} gave a
specific construction, for a Lie group $G$, of a space homotopy
equivalent to $BG$.  The basic fact we will require concerning
classifying spaces is the following:

\begin{Thm} If $H$ is a closed subgroup of the Lie group $G$, then there exists a natural map
$BH\rightarrow BG$ between the models constructed by Milnor; furthermore this map is a fiber
bundle with fiber the homogenous space $G/H$.
\end{Thm}

\subsection{Okun's construction}

Okun established \cite[Theorem 5.1]{O1} the following nice result:

\begin{Thm}\label{Theorem2.2} Let $X=\G \bs G/K$ and $X_u=G_u/K$ be dual symmetric spaces.
Then there exists a finite sheeted cover $\bar X$ of $X$ (ie a subgroup $\bar \G$ of
finite index in $\G$, $\bar X=\bar \G \bs G/K$), and a tangential map $k\co\bar X
\rightarrow X_u$.
\end{Thm}

This was subsequently used by Okun to exhibit exotic smooth
structures on certain compact locally symmetric spaces of
noncompact type \cite{O2}, and by Aravinda--Farrell in their
construction of exotic smooth structures on certain quaternionic
hyperbolic manifolds supporting metrics of strict negative curvature
\cite{AF}.  More recently, this was used by Lafont--Roy \cite{LaR} to
give an elementary proof of the Hirzebruch proportionality principle
for Pontrjagin numbers, as well as (non)vanishing results for
Pontrjagin numbers of the Gromov--Thurston examples of manifolds with
negative sectional curvature.

Since it will be relevant to our proof of the main theorem, we
briefly recall the construction of the finite cover that appears in
Okun's argument for \fullref{Theorem2.2}.  Starting from the canonical
principle fiber bundle
$$\G \bs G\rightarrow \G \bs G/K=X$$
with structure group $K$ over the base space $X$, we can extend the
structure group to the group $G$, yielding the flat principle
bundle:
$$\G \bs G\times_K G=G/K \times_\G G \longrightarrow \G \bs G/K=X$$
Further extending the structure group to $G_\mC$ yields a flat
bundle with a complex linear algebraic structure group.  A result of
Deligne and Sullivan \cite{DS} implies that there is a finite cover
$\bar X$ of $X$ where the pullback bundle is trivial; since $G_u$ is
the maximal compact in $G_\mC$, the bundle obtained by extending the
structure group from $K$ to $G_u$ is trivial as well. In terms of
the classifying spaces, this yields the commutative diagram:
$$\xymatrix{ & G_u/K \ar[d]\\
\bar X \ar[r] \ar@{.>}[ur] \ar[dr]_{\simeq 0}& BK \ar[d]\\
 & BG_u}$$
Upon homotoping the bottom diagonal map to a point, one obtains that
the image of the horizontal map lies in the fiber above a point,
ie inside $G_u/K$, yielding the dotted diagonal map in the above
diagram.  Okun then proceeds to show that the map to the fiber is
the desired tangential map (since the pair of maps to $BK$ classify
the respective canonical $K$--bundles on $\bar X$ and $G_u/K$, and
the canonical $K$--bundles determine the respective tangent bundles).

\subsection{Matsushima's map}

Matsushima \cite{Mat} constructed a map on cohomology
$j^*\co H^*(G_u/K; \mR)\rightarrow H^*(X; \mR)$ whenever $X$ is a
compact locally symmetric space modelled on $G/K$.  We will require the
following fact concerning the Matsushima map:

\begin{Thm}[Matsushima \cite{Mat}] The map $j^*$ is always injective.
Furthermore, there exists a constant $m(\frakg)$ (called the Matsushima
constant) depending solely on the Lie algebra $\frakg$ of the Lie group
$G$, with the property that the Matushima map $j^*$ is a surjection
in cohomology up to the dimension $m(\frakg)$.
\end{Thm}

The specific value of the Matsushima constant for
the locally symmetric spaces that are K\"ahler can be found in \cite{Mat}.
We also point out the following result of Okun \cite[Theorem
6.4]{O1}:

\begin{Thm}\label{Theorem2.4} Let $X=\G \bs G/K$ be a compact locally symmetric space,
and $\bar X$, $t\co \bar X \rightarrow G_u/K$ the finite cover and
tangential map constructed in \fullref{Theorem2.2}.  If the groups $G_u$ and
$K$ have equal rank, then the induced map $t^*$ on cohomology
coincides with Matsushima's map $j^*$.
\end{Thm}

\section{Detecting homologically essential submanifolds}

In this section, we provide a proof of \fullref{Theorem1.1}, which gives a criterion for establishing
when a totally geodesic submanifold $Y\subset X$ in a locally symmetric space of noncompact
type, is homologically nontrivial.

\begin{proof}[Proof of \fullref{Theorem1.1}]
In order to establish the theorem, we make use of differential
forms.  If a group $H$ acts on a smooth manifold $M$, we let $\gO
^H(M)$ denote the complex of $H$--invariant differential forms on
$M$. Let $X=\Gamma \bs G/K$, $Y=\gL \bs G^\p/K^\p$ be the pair of
compact locally symmetric spaces, and $X_u=G_u/K$, $Y_u=G_u^\p/K^\p$
be the corresponding dual spaces.  We now consider the following
four complexes of differential forms: (1) $\gO ^{G}(G/K)$, (2) $\gO
^{G^\p}(G^\p/K^\p)$, (3) $\gO ^\G(G/K)$ and (4) $\gO
^\gL(G^\p/K^\p)$.

We now observe that the cohomology of the first two complexes can be identified with the
cohomology of
$X_u$, $Y_u$ respectively.  Indeed, we have the sequence of natural identifications:
$$\gO ^{G}(G/K)=H^*(\frakg,\ft)=H^*(\frakg_u, \ft)=\gO ^{G_u}(G_u/K)$$
The first and third equalities come from the identification of the complex of harmonic
forms with the relative Lie algebra cohomology.  The second equality comes via the dual Cartan
decompositions: $\frakg=\ft \oplus \fp$ and $\frakg_u=\ft \oplus i\fp$.  Since $X_u=G_u/K$ is a compact closed
manifold, and $\gO ^{G_u}(G_u/K)$ is the
complex of harmonic forms on $X_u$, Hodge theory tells us that the cohomology of the
complex $\gO ^{G_u}(G_u/K)$ is just the cohomology of $X_u$.
The corresponding analysis holds for $\gO ^{G^\p}(G^\p/K^\p)$.

Next we note that the cohomology of the last two complexes can be
identified with the cohomology of $X$, $Y$ respectively.  This just
comes from the fact that the projection $G/K\rightarrow \G \bs
G/K=X$ induces the isomorphism of complexes $\gO ^\G(G/K)=\gO(X)$,
and similarly for $Y$.

Now observe that the four complexes fit into a commutative diagram
of chain complexes:
$$\xymatrix{\gO ^\gL(G^\p/K^\p) & \gO ^\G(G/K)\ar[l]_{~\phi}\\
\gO ^{G^\p}(G^\p/K^\p) \ar[u]^{j_Y} & \gO ^{G}(G/K)\ar[l]_{~\psi}
\ar[u]_{j_X}}$$ Let us briefly comment on the maps in the diagram.
The vertical maps are obtained from the fact that $\G\leq G$, so
that any $G$--invariant form can be viewed as a $\G$--invariant form,
and similarly for $\gL \leq G^\p$.

For the horizontal maps, one observes that $G^\p/K^\p
\hookrightarrow G/K$ is an embedding, hence any form on $G/K$
restricts to a form on $G^\p/K^\p$.  We also have the inclusion
$\gL\leq \G$, and hence the restriction of a $\G$--invariant form on
$G/K$ yields a $\gL$--invariant form on $G^\p/K^\p$.  This is the
horizontal map in the top row.  One obtains the horizontal map in
the bottom row similarly.

Now passing to the homology of the chain complexes, and using the
identifications discussed above, we obtain a commutative diagram in
dimension $m=\dim(Y)=\dim(G_u^\p/K^\p)$:
$$\xymatrix{\mR \simeq H^m(Y; \mR) & & H^m(X;\mR) \ar[ll]_{\phi^*}\\
\mR \simeq H^m(G_u^\p/K^\p;\mR) \ar[u]^{j_Y^*} & & H^m(G_u/K;\mR)
\ar[ll]_{\psi^*} \ar[u]_{j_X^*}}$$ Note that the two vertical maps defined
here are precisely the Matsushima maps for the respective locally
symmetric spaces.  Since Matsushima's map is always injective, and
the cohomology of $H^m(G_u^\p/K^\p;\mR)$ and $H^m(Y; \mR)$ are both
one-dimensional, we obtain that $j_Y^*$ is an isomorphism.  Likewise
$j_X^*$ is always injective, and if $m\leq m(\frakg)$ then $j_X^*$ is also
surjective (and hence $j_X^*$ is an isomorphism as well).  This implies the 
following two facts:
\begin{itemize}
\item If $\phi^*$ is identically zero, then $\psi^*$ is identically zero.
\item If furthermore $m\leq m(\frakg)$, then both vertical maps are isomorphisms,
and we have that $\psi^*$ is identically zero if and only if $\phi^*$ is identically zero.
\end{itemize}
Now observe that
both of the horizontal maps coincide with the maps induced on
cohomology by the respective inclusions $Y\hookrightarrow X$ and
$G_u^\p/K^\p \hookrightarrow G_u/K$; indeed the maps are obtained by
restricting the forms defined on the ambient manifold to the
appropriate submanifold.  In particular, the map $\psi^*$ coincides with
the map $\rho$ that appears in the statement of our theorem.  On the other hand, from 
the Kronecker pairing, the map $\phi^*$ is nonzero precisely when
$[Y^m]\neq 0\in H_m(X;\mR)$.  Combining these observations with the
two facts in the previous paragraph completes the proof of \fullref{Theorem1.1}.
\end{proof}

\begin{Rmk}
(1)\qua The Matsushima map is only defined on the real cohomology
(since it passes through differential forms), and as a result, {\it
cannot} be used to obtain any information on torsion elements in
$H^k(X^n;\mZ)$.

(2)\qua We remark that the proof of \fullref{Theorem1.1} applies equally well to 
lower-dimen\-sional cohomology (using the fact that Matsushima's map is injective
in all dimensions), and gives the following lower-dimensional criterion.  Assume
that
the map $H^k(X^n_u; \mR) \rightarrow H^k(Y^m_u;\mR)$ has image containing a
nonzero class $\alpha$, and let $i(\alpha)\in H^k(Y^m; \mR)$ be the nonzero image class
under the Matsushima map.  Then the homology class $\beta\in H_k(Y^m;\mR)$ dual
(under the Kronecker pairing) to $i(\alpha)$ has nonzero image in $H_k(X^n;\mR)$ under
the map induced by the inclusion $Y^m\hookrightarrow X^n$.
\end{Rmk}

\section{Pairs of tangential maps}
In this section, we proceed to give a proof of \fullref{Theorem1.2},
establishing the existence of pairs of tangential maps from the pair
$(\bar X,\bar Y)$ to the pair $(X_u,Y_u)$.

\begin{proof}[Proof of \fullref{Theorem1.2}]
We start out by applying
\fullref{Theorem2.2}, which gives us a finite cover $\bar X$ of $X$ with the
property that the natural composite map $\bar X\rightarrow
BK\rightarrow BG_u$ is homotopic to a point. Note that this map
classifies the principle $G_u$ bundle over $\bar X$.

Now let $\bar Y\hookrightarrow \bar X$ be a connected lift of the
totally geodesic subspace $Y\hookrightarrow X$.  Observe that, by
naturality, we have a commutative diagram:
$$
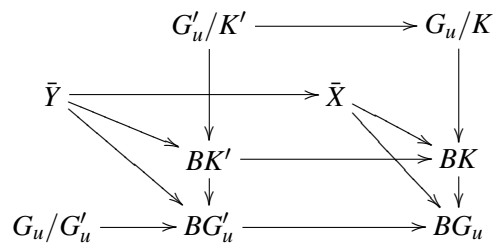
\xymatrix@R=10pt{ & G_u^\p/K^\p \ar[rr] \ar[dd] & & G_u/K \ar[dd]\\
\bar Y \ar[rr] \ar[dr] \ar[ddr] & & \bar X \ar[dr] \ar[ddr] & \\
& BK^\p \ar[rr] \ar[d] & & BK \ar[d] \\
G_u/G_u^\p \ar[r] & BG_u^\p \ar[rr] & & BG_u}$$ 
By Okun's result,
the composite map $\bar X\rightarrow BG_u$ is homotopic to a point
via a homotopy $H\co\bar X \times I\rightarrow BG_u$.  We would like
to establish the existence of a homotopy $F\co\bar Y\times I
\rightarrow BG_u^\p$ with the property that the following diagram
commutes:
$$\xymatrix{\bar Y\times I \ar[r]^{i\times \mathrm{Id}} \ar[d]_F & \bar X\times I \ar[d]^H\\
BG_u^\p \ar[r] & BG_u}$$
Indeed, if we had the existence of such a compatible pair of
homotopies, then one can easily complete the argument: since each of
the vertical columns in the diagram are fiber bundles, we see that
after applying the pair of compatible homotopies, the images of
$(\bar X,\bar Y)$ lies in the pair of fibers $(G_u/K, G_u^\p/K^\p)$.
This yields a pair of compatible lifts, yielding a commutative
diagram of the form:
$$\xymatrix{ & G_u^\p/K^\p \ar[rr] \ar[dd] & & G_u/K \ar[dd]\\
\bar Y \ar[ru] \ar[rr] \ar[dr]  & & \bar X \ar[dr] \ar[ru] & \\
 & BK^\p \ar[rr] & & BK } $$
Since the pair of maps to $BK^\p$ (respectively $BK$) classify the
canonical $K^\p$--bundle structures on $\bar Y$, $G_u^\p/K^\p$
(respectively the canonical $K$--bundle structure on $\bar X$,
$G_u/K$), and since these bundles canonically determine the tangent
bundles of these spaces \cite[Lemma 2.3]{O1}, commutativity of the
diagram immediately
gives us tangentiality of the maps $\bar Y \rightarrow G_u^\p/K^\p$
(respectively, of the map $\bar X\rightarrow G_u/K$).

In order to show the existence of the compatible homotopy $F\co\bar
Y\times I \rightarrow BG_u^\p$, we start by observing that the
bottom row of the commutative diagram is in fact a fibration
$$ G_u/G_u^\p \rightarrow BG_u^\p \rightarrow BG_u.$$
Since $\bar Y$ is embedded in $\bar X$, we see that the homotopy $H$
induces by restriction a homotopy $H\co\bar Y \times I \rightarrow
BG_u$.  Since the bottom row is a fibration, we may lift this
homotopy to a homotopy $\tilde H\co \bar Y\times I\rightarrow
BG_u^\p$, with the property that $\tilde H_0$ coincides with the map
$\bar Y \rightarrow BG_u^\p$ which classifies the canonical
principle $G_u^\p$ bundle over $\bar Y$.  Unfortunately, we do not
know, a priori, that the time one map $\tilde H_1$ maps $\bar Y$ to
a point in $BG_u^\p$.  Indeed, we merely know that $\tilde H_1(\bar
Y)$ lies in the preimage of a point in $BG_u$, ie in the fiber
$G_u/G_u^\p$.  Our next goal is to establish that the map $\tilde
H_1\co \bar Y \rightarrow G_u/G_u^\p$ is nullhomotopic.  If this were
the case, we could concatenate the homotopy $\tilde H$ taking $\bar
Y$ into the fiber $G_u/G_u^\p$ with a homotopy contracting $\tilde
H_1\co\bar Y \rightarrow G_u/G_u^\p$ to a point {\it within the
fiber}.  This would yield the desired homotopy $F$.
\medskip

In order to establish that $\tilde H_1\co \bar Y \rightarrow
G_u/G_u^\p$ is nullhomotopic, we merely note that we have the
fibration
$$G_u^\p\rightarrow G_u\rightarrow G_u/G_u^\p.$$
From the corresponding long exact sequence in homotopy groups, and
using the fact that the inclusion $G_u^\p\hookrightarrow G_u$
induces a $\pi_i$--isomorphism for $i< m$ and a surjection on
$\pi_m$, we immediately obtain that $\pi_i(G_u/G_u^\p)\cong 0$ for
$i\leq m$. Since the dimension of the manifold $\bar Y$ is $m$, we
can now conclude that the map $\tilde H_1$ is nullhomotopic.
Indeed, taking a cellular decomposition of $\bar Y$ with a single
$0$--cell, one can recursively contract the image of the $i$--skeleton
to the image of the $0$--cell: the obstruction to doing so lies in
$\pi_i(G_u/G_u^\p)$, which we know vanishes.  This yields that
$\tilde H_1$ is nullhomotopic, which by our earlier discussion,
implies the existence of a tangential map of pairs $(\bar X,\bar Y)
\rightarrow (X_u,Y_u)$. Finally, to conclude we merely point out the
Okun has shown (see \fullref{Theorem2.4}) that in the case where the rank of
$G_u$ equals the rank of $K$, the tangential map he constructed
induces the Matsushima map on cohomology.  Our construction
restricts to Okun's construction on both $X$ and $Y$, and from the
hypothesis on the ranks, so we conclude that the tangential map of
pairs induces the Matsushima map on the cohomology of each of the
two spaces. This concludes the proof of \fullref{Theorem1.2}.
\end{proof}
\medskip

\begin{Rmk}
We observe that the argument given above, for the case of a pair $(X^n, Y^m)$,
can readily be adapted to deal with any descending chain of totally geodesic
submanifolds.  More precisely, assume that we have a series of totally
geodesic embeddings $X^n=Y_k \supset \cdots \supset Y_2\supset Y_1$, with
the property that each $Y_j$ is a closed locally symmetric space of noncompact type.
Further assume that, if $(Y_j)_u = (G_j)_u/K_j$ denotes the compact duals,
the maps $(G_j)_u \hookrightarrow (G_{j+1})_u$
induced by the inclusions $Y_j\hookrightarrow Y_{j+1}$ are $\pi_i$ isomorphisms
for $i< \dim(Y_j)$ and a surjection on $\pi_i$ ($i=\dim(Y_j)$).  Then there exists
a finite cover $\bar X^n=\bar Y_k$ of $X^n$, and connected lifts $\bar Y_j$ of $Y_j$,
having the property that:
\begin{itemize}
\item we have containments $\bar Y_j\subset \bar Y_{j+1}$, and
\item there exists a map $(\bar Y_k, \ldots, \bar Y_1)\rightarrow \big((Y_k)_u,
\ldots ,(Y_1)_u\big)$ which restricts to a tangential map from each $\bar Y_j$ to the
corresponding $(Y_j)_u$.
\end{itemize}
This is shown by induction on the length of a descending chain.  We leave the details
to the interested reader.
\end{Rmk}

We now proceed to show \fullref{Corollary1.4}, that is to say, that in the
case where $X^n$ is real hyperbolic, and $Y^m\hookrightarrow X^n$ is
totally geodesic, there exists a finite cover $\bar X$ of $X^n$ and
a connected lift $\bar Y$ of $Y^m$, with the property that the
normal bundle of $\bar Y$ in $\bar X$ is trivial.

\begin{proof}[Proof of \fullref{Corollary1.4}] We first observe that, provided one could verify the
hypotheses of \fullref{Theorem1.2} for the pair $(X^n,Y^m)$, the corollary
would immediately follow.  Note that in this case, the dual spaces
$X_u$ and $Y_u$ are spheres of dimension $n$ and $m$ respectively.
This implies that the totally geodesic embedding $Y_u\hookrightarrow
X_u$ is in fact a totally geodesic embedding $S^m\hookrightarrow
S^n$, forcing the normal bundle to $Y_u$ in $X_u$ to be trivial. But
now \fullref{Corollary1.3} to the \fullref{Theorem1.2} immediately yields \fullref{Corollary1.4}.

So we are left with establishing the hypotheses of \fullref{Theorem1.2} for
the pair $(X^n,Y^m)$.  We observe that in this situation we have
the groups $G_u \cong SO(n+1)$, and $G_u^\p \cong SO(m+1)$.
Furthermore, there is essentially a unique totally geodesic
embedding $S^m\hookrightarrow S^n$, hence we may assume that the
embedding $G_u^\p\hookrightarrow G_u$ is the canonical one.  But now
we have the classical facts that (1) the embeddings
$SO(m+1)\hookrightarrow SO(n+1)$ induce isomorphisms on $\pi_i$ for
$i< m$ and (2) that the embedding induces a surjection
$\pi_m(SO(m+1))\rightarrow \pi_m(SO(n+1))$.  Indeed, this is
precisely the range of dimensions where the homotopy groups
stabilize \cite{Mil2}.  This completes the verification of the
hypotheses, and hence the proof of \fullref{Corollary1.4}.
\end{proof}

We now proceed to give an example of an inclusion
$Y^m\hookrightarrow X^n$ satisfying the hypotheses of our theorem.
Our spaces will be modelled on complex hyperbolic spaces, namely we
have:
\begin{align*}Y^{2m}&= \Lambda\backslash \mC\mH ^m = \Lambda\backslash \SU(m,1) / S(U(m)\times U(1))\\
X^{2n}& = \Gamma\backslash \mC\mH ^n = \Gamma\backslash \SU(n,1) / S(U(n)\times U(1))\end{align*}
To construct such pairs, one starts with the standard inclusion of $\SU(m,1)\hookrightarrow \SU(n,1)$,
which induces a totally geodesic embedding $\mC\mH ^m\hookrightarrow \mC\mH ^n$.  One can
now construct explicitly (by arguments similar to those in Gromov and Piatetski-Shapiro \cite{GP}) an arithmetic uniform lattice
$\Lambda \leq \SU(m,1)$ having an extension
to an arithmetic uniform lattice $\Gamma \leq \SU(n,1)$.  Quotienting out by these lattices gives the
desired pair.

Let us now consider these examples in view of our \fullref{Theorem1.2}.  First of all, we have that the
respective complexifications are $G^\p_\mC=\SL(m+1, \mC)$ and $G_\mC=\SL(n+1,\mC)$, with
the natural embedding
$$G^\p_\mC =\SL(m+1,\mC) \hookrightarrow \SL(n+1,\mC)=G_\mC.$$
Looking at the respective maximal compacts, we see that $G^\p_u=\SU(m+1)$, $G_u=\SU(n+1)$, and the inclusion is
again the natural embeddings
$$G^\p_u = \SU(m+1) \hookrightarrow \SU(n+1) = G_u.$$
Hence the homotopy condition in our theorem boils down to asking
whether the natural embedding $\SU(m+1)\hookrightarrow \SU(n+1)$
induces isomorphisms on the homotopy groups $\pi_i$, where $i\leq
\dim(Y^{2m})=2m$.  But it is a classical fact that the natural
embedding induces isomorphisms in all dimensions $i<2(m+1)=2m+2$,
since this falls within the stable range for the homotopy groups
(and indeed, one could use complex Bott periodicity to compute the
exact value of these homotopy groups \cite{Mil2}).  Finally, we observe
that in this context, the dual spaces are complex projective spaces, and the
embedding of dual spaces is the standard embedding $\mC`P^m\hookrightarrow \mC`P^n$.
It is well known that for the standard embedding, we have that the induced map on
cohomology $H^*(\mC`P^n) \rightarrow H^*(\mC`P^m)$ is surjective on cohomology.
Our \fullref{Theorem1.1} now tells us that $Y^{2m}\hookrightarrow X^{2n}$ is homologically
nontrivial.  Furthermore, we note that for these manifolds, $\rk(G_u)=\rk(K)$ and
$\rk(G_u^\p)=\rk(K^\p)$, and hence \fullref{Theorem1.2} tells us that
the cohomological map from the proof of \fullref{Theorem1.1}
can be (rationally) realized via a tangential map of pairs.

\section{Concluding remarks}

We conclude this paper with a few comments and questions.  First of all, in
view of our \fullref{Theorem1.1}, it is reasonable to ask for the converse:

\medskip
{\bf Question}\qua  Given an element $\alpha \in H_m(X^n
;\mR)$, is there an $m$--dimensional totally geodesic submanifold
$Y^m$ with $[Y^m]=\alpha$?
\medskip

A cautionary example for the previous question is provided by the following:

\begin{Prop}
Let $X$ be a compact hyperbolic $3$--manifold that fibers over $S^1$, with
fiber a surface $F$ of genus $\geq 2$.  Then the homology class
represented by $[F]\in H_2(X; \mZ)$ cannot be represented by a
totally geodesic submanifold.
\end{Prop}

\begin{proof}
Assume that there were such a totally geodesic submanifold $Y
\subset X$, and observe that since $Y$ is totally geodesic, we have
an embedding $\pi_1(Y)\hookrightarrow \pi_1(X)$.  Furthermore, since
$X$ fibers over $S^1$ with fiber $F$, we also have a short exact
sequence:   
$$ 0\rightarrow \pi_1(F)\rightarrow \pi_1(X)\rightarrow \mZ \rightarrow 0 $$
Our goal is to show that $\pi_1(Y)\subset \pi_1(F)$.  Indeed, if we
could establish this containment, one could then argue as follows:
since $Y$ is a compact surface, covering space theory implies that
$\pi_1(Y)\subset \pi_1(F)$ is a finite index subgroup.  Now pick a
point $x$ in the universal cover $\tilde X\cong \mH^3$, and consider
the subset $\Lambda_Y \subset \partial ^\infty \mH^3=S^2$ obtained
by taking the closure of the $\pi_1(Y)$--orbit of $x$.  Since $Y$ is
assumed to be totally geodesic, the subset $\Lambda\subset S^2$ is a
tamely embedded $S^1$ (identified with the boundary at infinity of a
suitably chosen totally geodesic lift $\tilde Y \cong \mH^2$ of
$Y$). On the other hand, since $\pi_1(Y)$ has finite index in
$\pi_1(F)$, we have that $\Lambda_Y$ must coincide with $\Lambda_F$,
the closure of the $\pi_1(F)$--orbit of $x$.  But the latter, by a
well-known result of Cannon--Thurston is known to be the entire
boundary at infinity (see for instance Mitra \cite{Mit}).

So we are left with establishing that $\pi_1(Y)\subset \pi_1(F)$. In
order to see this, let us consider the cohomology class $\alpha_F\in
H^1(X; \mZ)$ which is Poincar\'e dual to the class $[F]\in H_2(X;
\mZ)$. Now recall that the evaluation of the cohomology class
$\alpha_F$ on an element $[\gamma]\in H_1(X;\mZ)$ can be interpreted
geometrically as the intersection number of the representing curve
$\gamma$ with the surface $F$.  Furthermore, we have that the group
$H_1(X;\mZ)$ is generated by the image of $H_1(F;\mZ)$, under the
inclusion $F\hookrightarrow X$, along with an element $[\eta]\in
H_1(X; \mZ)$ mapping to $[S^1]\in H_1(S^1;\mZ)$.  Here $\eta$ is
chosen to be a closed loop in $M$ with the property that $\eta$ maps
homeomorphically to the base $S^1$ (preserving orientations) under
the projection map.  This gives us the following two facts:
\begin{itemize}
\item The class $\alpha_F$ evaluates to $1$ on the element $[\eta]$, since
$F\cap \eta$ is a single transverse point.
\item The class $\alpha_F$ evaluates to zero on the image of $H_1(F;\mZ)$ in
the group $H_1(X;\mZ)$.  This follows from the fact that
the surface $F$ has trivial normal bundle in $X$, allowing any curve
in $F$ representing (the image of) an element in $H_1(F; \mZ)$ to be homotoped to a
curve disjoint from $F$.
\end{itemize}
Furthermore, since we are assuming that $[Y]=[F]\in H_2(X;\mZ)$, we know that
we have an identification of Poincar\'e duals $\alpha_Y=\alpha_F$.

Now let us assume that $\pi_1(Y)$ is {\it not\/} contained in $\pi_1(F)$, and
observe that this implies that there exists a closed loop $\gamma\subset Y$ having
the property that under the composition $Y\hookrightarrow X\rightarrow S^1$,
the class $[\gamma]\in H_1(Y;\mZ)$ maps to $k\cdot [S^1]\in H_1(S^1; \mZ)$ (and
$k\neq 0$).  We now proceed
to compute, in two different ways, the evaluation of the cohomology classes $\alpha_Y
=\alpha _F$ on a suitable multiple of the homology class $[\gamma]$.

Firstly, from the comments above, we can write $[\gamma]$ as the sum of
$k\cdot [\eta]$, along with an element $\beta$, where $\beta$ lies in the image
of  $H_1(F; \mZ)$.  By linearity of
the Kronecker pairing, along with the two facts from the previous paragraph, we
obtain:
$$\alpha _F([\gamma])=\alpha_F(\beta) + k\alpha_F([\eta]) = k \neq 0$$
Secondly, observe that $Y$ is assumed to be embedded in $X$, and
represents the nonzero homology class $[Y]=[F]\in H_2(X;Z)$. This
implies that $Y$ must be orientable, and hence has trivial normal
bundle in $X$.  In particular, the curve $\gamma \subset Y$ can be
homotoped (within $X$) to have image disjoint from $Y$.  Since the
integer $\alpha _Y([\gamma])$ can be computed geometrically as the
intersection number of the curve $\gamma$ with $Y$, we conclude that
$\alpha_Y([\gamma])=0$.

Combining the two observations above, we see that
$0=\alpha_Y([\gamma])= \alpha_F([\gamma])\neq 0$, giving us the
desired contradiction. This completes the proof of the proposition.
\end{proof}

We observe that Thom \cite{T} has shown that in dimensions $0\leq
k\leq 6$ and $n-2\leq k\leq n$, every integral homology class can be
represented by an immersed submanifold.  In general however, there
can exist homology classes which are {\it not\/} representable by
submanifolds (see for instance Bohr, Hanke and
Kotschick \cite{BHK}).  The question above asks for a more stringent
condition, namely that the immersed submanifold in question be
totally geodesic.  We believe that the weaker question is also of
some interest, namely:

\medskip
{\bf Question}\qua  Find an example $X^n$ of a compact locally
symmetric space of noncompact type, and a homology class in some
$H_k(X^n; \mZ)$ which {\it cannot\/} be represented by an immersed
submanifold.
\medskip

Now our original motivation for looking at totally geodesic submanifolds
inside locally symmetric spaces was the desire to exhibit lower-dimensional
bounded cohomology classes.  In \cite{LaS}, the authors showed that for
the fundamental group $\Gamma$ of
a compact locally symmetric space of noncompact type $X^n$, the comparison
map from bounded cohomology to ordinary cohomology:
$$H^*_b(\Gamma)\rightarrow H^*(\Gamma)$$
is {\it surjective\/} in dimension $n$.  The proof actually passed
through the dual formulation, and showed that the $L^1$
(pseudo)-norm of the fundamental class $[X^n]\in H_n(X^n; \mR)$ is
nonzero.  Now given a totally geodesic embedding $Y\hookrightarrow
X$ of the type considered in this paper, it is tempting to guess
that the homology class $[Y]$ also has nonzero $L^1$ (pseudo)-norm.
Of course, this naive guess fails, since one can find examples where
$[Y]=0\in H_m(X^n;\mR)$.  The problem is that despite the fact that
the {\it intrinsic\/} $L^1$ (pseudo)-norm of $[Y]$ is nonzero, the
{\it extrinsic\/} $L^1$ (pseudo)-norm of $[Y]$ is zero.  In other
words, one can represent the fundamental class of $Y$ {\it more
efficiently} by using simplices that actually {\it do not lie in
$Y$} (despite the fact that $Y$ is totally geodesic).  The authors were unable to answer the following:

\medskip
{\bf Question}\qua Assume that $Y$ and $X$ are compact locally
symmetric spaces of noncompact type, that $Y\subset X$ is a totally
geodesic embedding, and that $Y$ is orientable with  $[Y]\neq 0 \in
H_m(X^n; \mR)$.  Does it follow that the dual cohomology class
$\beta\in H^m(X^n; \mR)$ (via the Kronecker pairing) has a bounded
representative?
\medskip

Now one situation in which nonvanishing of the $L^1$ (pseudo)-norm would
be preserved is the case where $Y\hookrightarrow X$ is actually a retract of $X$.
Hence one can ask the following:

\medskip
{\bf Question}\qua If $Y\subset X$ is a compact totally
geodesic proper submanifold inside a locally symmetric space of
noncompact type, when is $Y$ a retract of $X$?

\begin{Rmk}
(1)\qua  In the case where $X$ is a (nonproduct) higher rank locally
symmetric space of noncompact type, one cannot find a proper
totally geodesic submanifold $Y\subset X$ which is a retract of $X$.
Indeed, if there were such a submanifold, then the morphism
$\rho\co\pi_1(X)\rightarrow \pi_1(Y)$ induced by the retraction would have
to be surjective.  By Margulis' normal subgroup theorem \cite{Mar},
this implies that either (1) $\ker(\rho)$ is finite, or (2)
the image $\pi_1(Y)$ is finite. Since $Y$ is locally symmetric
of noncompact type, $\pi_1(Y)$ cannot be finite, and hence we must
have finite $\ker(\rho)$. But $\ker(\rho)$ is a subgroup of the torsion-free
group $\pi_1(X)$, hence must be trivial.  This forces $\pi_1(X)\cong \pi_1(Y)$,
which contradicts the fact that the cohomological dimension of $\pi_1(X)$
is $\dim (X)$, while the cohomological dimension of $\pi_1(Y)$ is $\dim (Y)<
\dim (X)$.  This implies  that no such morphism exists, and hence no such
submanifold exists.  The authors thank C\,Leininger for pointing out this simple
argument.

(2)\qua In the case where $X$ has rank one, the question is more delicate.  Some
examples of such retracts can be found in a paper by Long and Reid \cite{LoR}.
We remark that in this case, the application to bounded cohomology is not too
interesting, as Mineyev \cite{Min} has already shown that the comparison map
in this situation is surjective in all dimensions $\geq 2$.
\end{Rmk}

Finally, we conclude this paper by pointing out that the Okun maps, 
while easy to define, are geometrically
very complicated.  We illustrate this with a brief comment on the 
{\it singularities\/} of the tangential maps
between locally symmetric spaces of noncompact type and their
nonnegatively curved dual spaces.  More precisely, for a smooth map
$f\co X\rightarrow X_u$, consider the subset $\Sing(f)\subset
X_u$ consisting of points $p\in X_u$ having the property that there
exists a point $q\in X$ satisfying $f(q)=p$, and $\ker(df(q))\neq 0$
(where $df\co T_qX\rightarrow T_pX_u$ is the differential of $f$ at the
point $q$).  We can now ask how complicated the set $\Sing(h)\subset X_u$ 
gets for $h$ a smooth map within the homotopy class of the Okun maps.

\begin{Prop}\label{Proposition5.2}
Let $X$ be a closed, locally symmetric space of noncompact type, 
$X_u$ the nonnegatively curved dual space, $t\co\bar
X\rightarrow \smash{X_u}$ the Okun map from a suitable finite
cover $\bar X$ of $X$, and $h$ an arbitrary smooth map in the
homotopy class of $t$.  Consider an arbitrary embedded compact
submanifold $N^k\subset X$, having the property that \textnormal{(1)} $[N^k]\neq 0\in
H^k(X_u; \mathbb Z)$, and \textnormal{(2)} $N^k$ is simply connected.
Then we have that $\Sing(h)\cap N^k\neq
\emptyset$.
\end{Prop}

\begin{proof}
Let $h$ be a smooth map homotopic to $t$, and assume that
$\Sing(h)\cap N^k=\emptyset$.  This implies that $dh$ has full
rank at every preimage point of $N^k\subset X_u$.  Choose a
connected component $S\subset \bar X$ of the set $h^{-1}(N^k)$,
and observe that the restriction of $h$ to $S$ provides a local
diffeomorphism (and hence a covering map) to $N^k$.  
Since $N^k$ is simply connected, $S$ is diffeomorphic to
$N^k$, and $h$ restricts to a diffeomorphism from $S\subset \bar X$
to $N^k\subset X_u$.

Next we observe that the homology class represented by $[S]\in
H_k(\bar X;\mZ)$ is nonzero, since this class has image, under $h$,
the homology class $[N_k]\neq 0\in H_k(X_u; \mZ)$.  But 
observe that $S\subset \bar
X$, being simply connected, supports a cellular decomposition with a 
single $0$--cell and {\it no $1$--cells\/}.  In particular,
the submanifold $S$ is homotopic to a point in the aspherical space $\bar X$, 
since one can recursively contract all the cells of dimension $\geq 2$ down to 
the $0$--cell.  This forces $[S]=0\in H_k(\bar X;\mZ)$, giving us the desired 
contradiction.
\end{proof}

We point out a simple example illustrating this last proposition.  If $X^{2n}$ is 
complex hyperbolic, then $X_u=\mC`P^n$, and one can take for $N^2=\mC`P^1$
the canonically embedded complex projective line.  Topologically, $N^2$ is 
diffeomorphic to $S^2$, hence is simply connected, while homologically we
have that $[N^2]$ is the generator for the cohomology group $H^2(\mC`P^n; 
\mZ) \cong \mZ$.  \fullref{Proposition5.2} now applies, and gives us that
{\it any\/} smooth map $h$ in the homotopy class of the Okun map must have 
singular set intersecting the canonical $\mC`P^1$.  A similar example appears
in the case of quaternionic hyperbolic manifolds, with the singular sets being
forced to intersect the canonical $\mO`P^1$ (diffeomorphic to $S^4$) 
inside the dual space $\mO`P^n$.

\bibliographystyle{gtart}
\bibliography{link}

\begin{thebibliography}{}
\providecommand\bibmarginpar{\leavevmode\marginpar}
\def\urlstyle#1{{\tt #1}}

\bibitem{AF}
\textbf{C\,S Aravinda}, \textbf{F\,T Farrell}, \emph{Exotic structures and
  quaternionic hyperbolic manifolds}, from: ``Algebraic groups and
  arithmetic'', Tata Inst. Fund. Res., Mumbai (2004)  507--524
  \xox{MR}{2094123}

\bibitem{BHK}
\textbf{C Bohr}, \textbf{B Hanke}, \textbf{D Kotschick},
  \href{http://dx.doi.org/10.1007/s002290200279} {\emph{Cycles, submanifolds,
  and structures on normal bundles}}, Manuscripta Math. 108 (2002) 483--494
  \xox{MR}{1923535}

\bibitem{DS}
\textbf{P Deligne}, \textbf{D Sullivan}, \emph{Fibr\'es vectoriels complexes
  \`a groupe structural discret}, C. R. Acad. Sci. Paris S\'er. A-B 281 (1975)
  Ai, A1081--A1083 \xox{MR}{0397729}

\bibitem{GP}
\textbf{M Gromov}, \textbf{I Piatetski-Shapiro}, \emph{Nonarithmetic groups in
  {L}obachevsky spaces}, Inst. Hautes \'Etudes Sci. Publ. Math.  (1988) 93--103
  \xox{MR}{932135}

\bibitem{H}
\textbf{D Husemoller}, \emph{Fibre bundles}, third edition, Graduate Texts in
  Mathematics 20, Springer, New York (1994) \xox{MR}{1249482}

\bibitem{LaR}
\textbf{J-F Lafont}, \textbf{R Roy}, \emph{A note on characteristic numbers of
  non-positively curved manifolds}, to appear in Expo. Math.

\bibitem{LaS}
\textbf{J-F Lafont}, \textbf{B Schmidt},
  \href{http://dx.doi.org/10.1007/s11511-006-0009-1} {\emph{Simplicial volume
  of closed locally symmetric spaces of non-compact type}}, Acta Math. 197
  (2006) 129--143

\bibitem{LoR}
\textbf{D Long}, \textbf{A Reid}, \emph{Subgroup separability and virtual
  retractions of groups}, to appear in Topology

\bibitem{Mar}
\textbf{G\,A Margulis}, \emph{Discrete subgroups of semisimple {L}ie groups},
  Ergebnisse series 17, Springer, Berlin (1991) \xox{MR}{1090825}

\bibitem{Mat}
\textbf{Y Matsushima}, \emph{On {B}etti numbers of compact, locally sysmmetric
  {R}iemannian manifolds}, Osaka Math. J. 14 (1962) 1--20 \xox{MR}{0141138}

\bibitem{Mil}
\textbf{J Milnor},
  \href{http://links.jstor.org/sici?sici=0003-486X(195605)2:63:3%3C430:COUBI%3%
E2.0.CO%3B2-5} {\emph{Construction of universal bundles. {II}}}, Ann. of Math.
  $(2)$ 63 (1956) 430--436 \xox{MR}{0077932}

\bibitem{Mil2}
\textbf{J Milnor}, \emph{Morse theory}, Based on lecture notes by M. Spivak and
  R. Wells. Annals of Mathematics Studies, No. 51, Princeton University Press,
  Princeton, N.J. (1963) \xox{MR}{0163331}

\bibitem{Min}
\textbf{I Mineyev}, \href{http://dx.doi.org/10.1007/PL00001686}
  {\emph{Straightening and bounded cohomology of hyperbolic groups}}, Geom.
  Funct. Anal. 11 (2001) 807--839 \xox{MR}{1866802}

\bibitem{Mit}
\textbf{M Mitra}, \emph{Cannon-{T}hurston maps for trees of hyperbolic metric
  spaces}, J. Differential Geom. 48 (1998) 135--164 \xox{MR}{1622603}

\bibitem{O1}
\textbf{B Okun}, \href{http://dx.doi.org/10.2140/agt.2001.1.709} {\emph{Nonzero
  degree tangential maps between dual symmetric spaces}}, Algebr. Geom. Topol.
  1 (2001) 709--718 \xox{MR}{1875614}

\bibitem{O2}
\textbf{B Okun}, \href{http://dx.doi.org/10.2140/agt.2002.2.381} {\emph{Exotic
  smooth structures on nonpositively curved symmetric spaces}}, Algebr. Geom.
  Topol. 2 (2002) 381--389 \xox{MR}{1917058}

\bibitem{T}
\textbf{R Thom}, \emph{Quelques propri\'et\'es globales des vari\'et\'es
  diff\'erentiables}, Comment. Math. Helv. 28 (1954) 17--86 \xox{MR}{0061823}

\end{thebibliography}

\end{document}